\documentclass[reqno]{jgokova}
\usepackage{epsfig,graphicx}

\newtheorem{Thm}{Theorem}
\newtheorem{Lem}[Thm]{Lemma}

\newtheorem{Rm}{Remark}
\newtheorem{Qt}{Question}

\def\C{{\mathbb C}}

\def\P{{\mathbb P}}

\newcommand{\ben}{\begin{enumerate}}
\newcommand{\een}{\end{enumerate}}
\newcommand{\be}{\begin{equation}}
\newcommand{\ee}{\end{equation}}
\newcommand{\bea}{\begin{eqnarray}}
\newcommand{\eea}{\end{eqnarray}}
\newcommand{\bc}{\begin{center}}
\newcommand{\ec}{\end{center}}

\newtheorem{thm}{Theorem}[section]
\newtheorem{cor}[thm]{Corollary}

\theoremstyle{definition}
\newtheorem{defn}[thm]{Definition}

\theoremstyle{remark}

\begin{document}

\setcounter{page}{1}

\title[Short title]{Cork twists and automorphisms of $3$-manifolds}
{\tiny {\it D\lowercase{edicated to} G\lowercase{eorge} F\lowercase{loyd} }}\\


\author[]{\small {Selman Akbulut} }

\thanks{Partially supported by NSF grant DMS 0905917}


\address{ G\"{o}kova Geometry Topology Institute,  Mu\u{g}la, T\"{u}rk\.{i}ye}
\email{akbulut.selman@gmail.com}

\begin{abstract}
Here we study two interesting smooth contractible manifolds, whose boundaries have non-trivial mapping class groups. The first one is a non-Stein contractible manifold, such that every self diffeomorphism of its boundary extends inside; implying that this manifold can not be a loose cork. The second example is a Stein contractible manifold which is a cork, with an interesting cork automorphism  $f:\partial W \to \partial W$.
By \cite{am} we know that any homotopy $4$-sphere is obtained gluing together two contractible Stein manifolds along their common boundaries by a diffeomorphism. We use the homotopy sphere $\Sigma = -W\smile_{f}W$ as a test case to investigate whether it is $S^4$? We show that
$\Sigma$ is a Gluck twisted $S^4$ twisted along a $2$-knot $S^{2}\hookrightarrow S^4$; then by  this we obtain a $3$-handle free handlebody description of $\Sigma$, and then show $\Sigma \approx S^4$.
\end{abstract}

\keywords{}

\maketitle

\vspace{-.04in}

\section{Introduction}

A {\it cork} is a pair $(W,f)$, where $W$ is a compact contractible Stein manifold, with an involution $f:\partial W\to \partial W$, which extends to a homeomorphism of $W$, but does not extend to a diffeomorphism of $W$.   We say $(W,f)$ is a cork of $M$,  if there is an imbedding 
$W\hookrightarrow M$  and cutting $W$ out of $M$ and re-gluing with $f $  produces an exotic copy $M'$ (\cite{ay}).

 $$M \mapsto M'=W\cup _{f}[M- W ]$$
 
The operation $M\mapsto M'$ is called {\it cork-twisting} $M$ along $W$. The first example of a cork appeared in \cite{a1}, then in \cite{m}, \cite{cfhs} it was proven that any exotic copy $M'$ of a closed simply connected $4$-manifold $M$ is obtained by twisting along a contractible manifold by an involution as above. Furthermore  in \cite{am}  it was shown that this contractible manifold can be taken to be a Stein manifold. If a Stein manifold $W$ has boundary $S^3$ then by EliashbergÕs theorem $W$ must be diffeomorhic to $B^4$, and hence it cannot be a cork.
 A cork without the ``Stein" condition is called a {\it loose-cork}. It is not known if loose-corks are corks (they have to contain corks by above).
 
 \begin{Qt}  \label{Q:loose}
 Is there any loose-cork with irreducible 
 boundary which can not be a cork?
 \end{Qt}

\begin{figure}[ht] 
 \begin{center}
 \includegraphics[width=.51\textwidth]{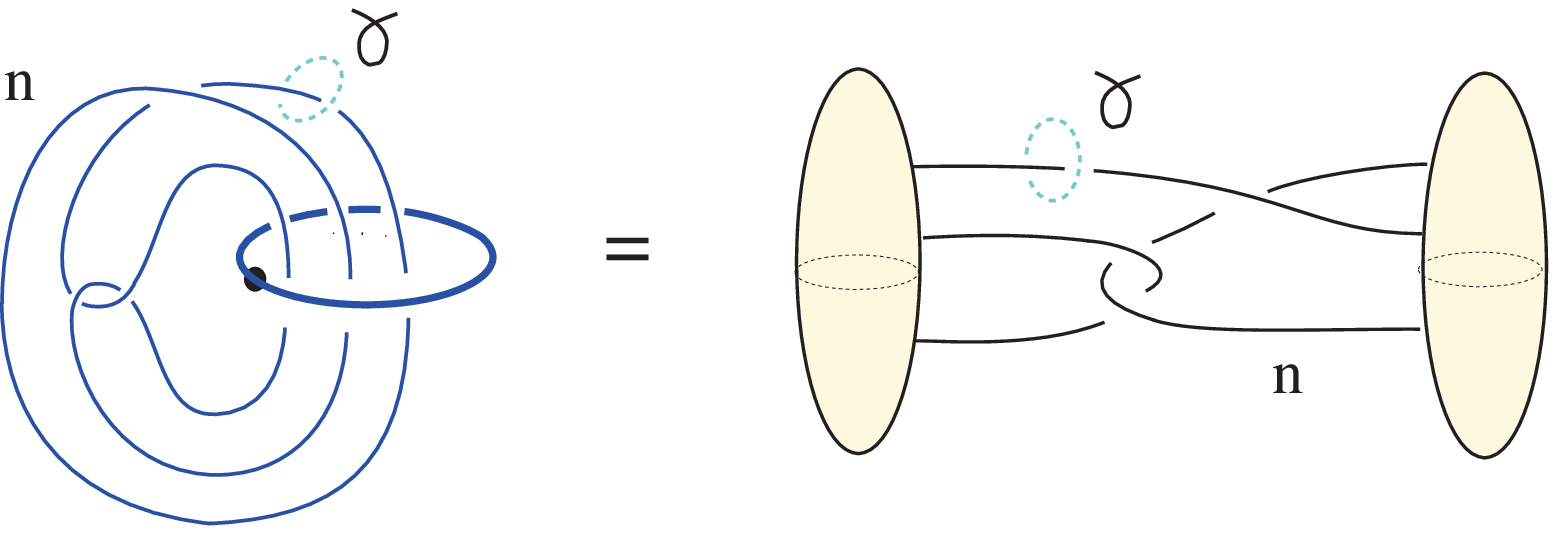}       
\caption{$W(0,n)$}     \label{c1} 
\end{center}
  \end{figure}

 The infinite order corks of  \cite{g}  and \cite{a2} could provide examples to this question. Recently, Mark and Tosun proved that the contractible manifold $W(0,2)$ can not be a Stein manifold \cite{mt}. In the notation of \cite{ak} and \cite{a3} this manifold is drawn in Figure \ref{c1} (n=2). Taking a connected sum of two copies gives an answer to Question~\ref{Q:loose} without the requirement of irreducibility (\cite[Corollary 1.8]{mt}). So the only other way a contractible Stein (or non-Stein) manifold $W$ fails to be a cork (or a loose-cork) is when all the self diffeomorphisms of $\partial W$ smoothly extend inside $W$. Here we test this on two interesting specific examples from the family $W(0,n)$, which were introduced in \cite{ak}.

From its Legendrian picture it is easy to check that when $n\leq 1$, $W(0,n)$ is a Stein manifold.  Moreover, $\partial W(0,1)$ and  $\partial W(0,2)$ can be identified as $+1$ and $-1$ surgeries of the Stevedore knot $K$, respectively. It is known that $K^{+1}$ is a hyperbolic manifold \cite{bw}, and $K^{-1}$ is the Brieskorn homology sphere   
$\Sigma(2,3,13)$ (\cite{ak}, \cite{a3}).
  \begin{figure}[ht]  
   \begin{center}
 \includegraphics[width=.25\textwidth]{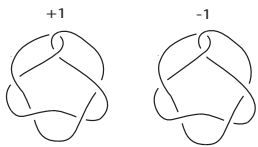}       
\caption{$\partial W(0,1)$ and $\partial W(0,2)$}      \label{c2} 
\end{center}
  \end{figure}
  
\begin{thm}\label{mainthm} 
The non-Stein manifold $W(0,2)$ can not be a loose-cork.
\end{thm}


We will construct a curious diffeomorphism $f: \partial W(0,1)\to \partial W(0,1)$ which is not obtained by the usual ``zero-dot exchange" process like many cork twisting maps 
. This leads us to the natural question whether the homotopy sphere  $\Sigma = -W(0,1)\smile_{f} W(0,1)$, formed by gluing two copies of $W(0,1)$  by $f$, is the standard smooth copy of $S^4$? Recall that, by \cite{am} we know every smooth homotopy sphere can be decomposed as a union of two contractible Stein manifolds glued along their boundaries (not necessarily by a cork twisting map). Therefore, it is only natural to seek counterexamples to $4$-dimensional Poincar\'{e} conjecture among such manifolds. Surprisingly, as the previous known examples  (\cite{a4},\cite{a5}), $\Sigma$ also decomposes as a Gluck twisted $S^4$, twisted along a knotted $S^2\hookrightarrow S^4$, (Figure~\ref{c52}).  By this, we cancel all $3$-handles of a handlebody of $\Sigma$. This process results a simple, $3$-handle-free handlebody
 description of $\Sigma$ (Figure~\ref{c58}), by which we prove 
 $\Sigma \approx S^{4}$. $f$ can be identified with the mapping class generators  $R$ or $S$ of $M(\partial W(0,1))$ in Lemma~\ref{P:mcg}.
 
 \begin{thm}\label{hs} 
$\Sigma$ is diffeomorphic to $S^{4}$.
\end{thm}

\section{Proofs}\label{S:extension}

Let $M(X)$ denote the mapping class group of a smooth manifod $X$.  It is orientation-preserving diffeomorphisms of $X$ modulo isotopy. 

\begin{figure}[ht]  \begin{center}
 \includegraphics[width=.4\textwidth]{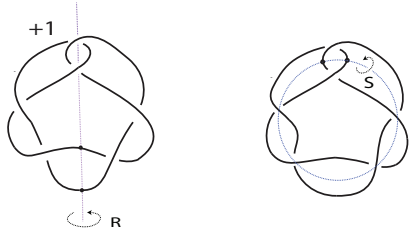}       
\caption{Self diffeomorphisms $R$ and $S$ of  $K^{+1}$}      \label{c3} 
\end{center}
\vspace{-.1in}
  \end{figure}

\begin{figure}[ht]  \begin{center}
 \includegraphics[width=.35\textwidth]{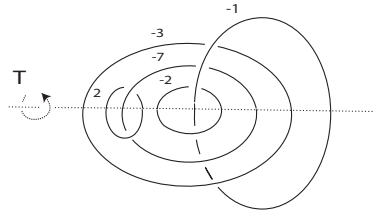}       
\caption{Self diffeomorphism $T$ of $K^{-1}$}      \label{c3prime} 
\end{center}
  \end{figure}

\begin{Lem}\label{P:mcg}
$M(K^{+1})$ is  $Z_{2}\oplus Z_{2}$, and is generated by the symmetries induced by the rotations $R$ and $S$ of the knot $K$, as indicated in Figure~\ref{c3}, and  $M(K^{-1})=Z_{2}$, which is generated by the symmetry $T$ of  Figure~\ref{c3prime}. In this figure we used another identification of $\Sigma(2,3,13)$ (from Exercise 12.3 of \cite{a3}), which is equivalent to this plumbing).
\end{Lem}


  \proof

$M(K^{-1})$ is straightforward to calculate, by using the fact that it is a Seifert fibered space.  From~\cite{ak}, we know that $K^{-1}$ is the Brieskorn sphere $\Sigma(2,3,13)$. This is a `small' Seifert fibered space, from which it follows that any element in $M(K^{-1})$ is isotopic to a fiber-preserving diffeomorphism~\cite{bo}. So, any orientation preserving self-diffeomorphism is isotopic to the identity, or to an involution that reverses the orientation of both base and fiber.  For the identification of $M(K^{+1})$ we refer reader to \cite{r}. \qed

\begin{figure}[ht]  \begin{center}
 \vspace{-.25in}
 \includegraphics[width=.5\textwidth]{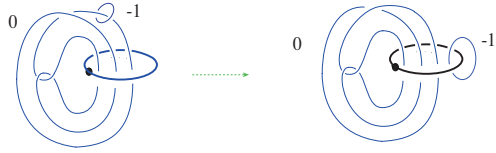}       
\caption{$W(0,1)\#\bar{C\P^2}$ and its absolutely exotic copy}      \label{c8} 
\end{center}
  \end{figure}

   \begin{Rm} Manifolds $W(0,n)$  have an interesting feature: Blowing them up $n$ times produces absolutely exotic manifold pairs with a cork inside (as in \cite{a6}). We will call contractible manifolds with this property ``{\it almost corks}". Figure~\ref{c8} demonstrates this process when $n=1$, iterating gives an absolutely exotic $W(0,n)\# n \C\P^{2}$ containing $W(0,0)$. \end{Rm}

\subsection{ Proof of Theorem 1}

\begin{figure}[ht]  \begin{center}
 \includegraphics[width=.63\textwidth]{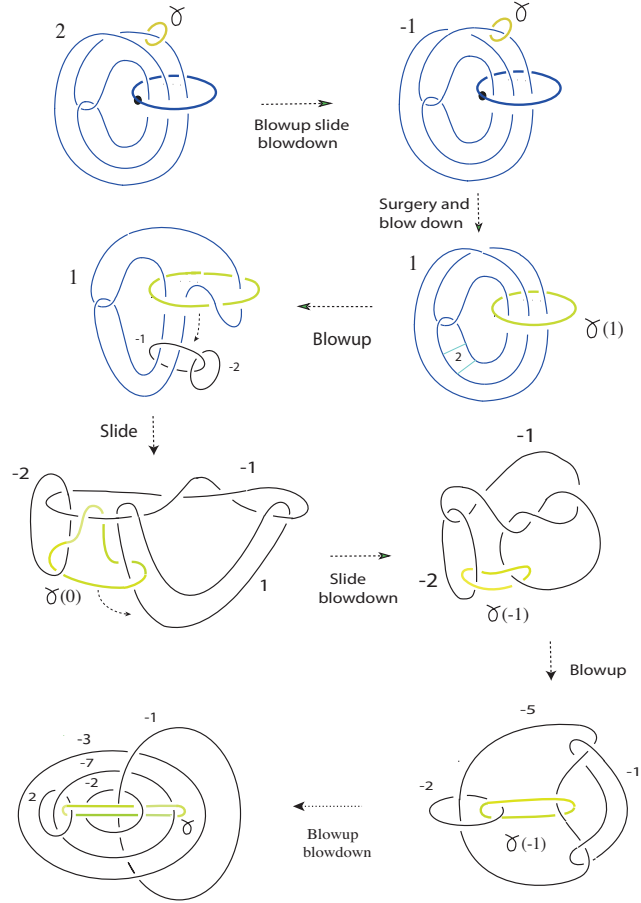}       
\caption{Tracing dual circle of $W(0,2)$ to $\Sigma(2,3,13)$}      \label{c6} 
\end{center}
  \end{figure}

\begin{figure}[ht]  \begin{center}
 \includegraphics[width=.58\textwidth]{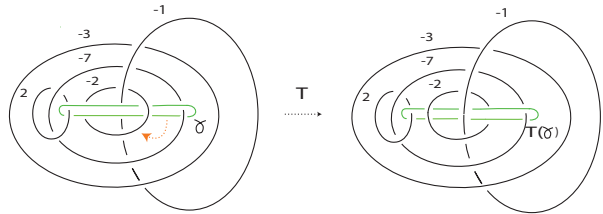}       
\caption{$T(\gamma)$ is isotopic to $\gamma$ }      \label{c3primeprime} 
\end{center}
  \end{figure}

\proof  We will show the generator $T$  of $M \big(\partial W(0,2)\big)$  extends to a self-diffeomorphism of $W(0,2)$. 
  It suffices to show that the diffeomorphism $T$ keeps the isotopy type of the dual loop of the $2$-handle of $W(0,2)$; then we can extend $T$ to a self diffeomorphism of $W(0,2)$ by carving (\cite{a3}). To see this we need to analyze the identification of $\partial W(0,2)$ with Figure~\ref {c3prime} closely. This is done in Figure~\ref{c6} (also see \cite{ak}). Now Figure~\ref{c3primeprime} shows that  the loop $T(\gamma)$ is isotopic to $\gamma $ (we can see this by sliding $T(\gamma)$ over the -2 framed $2$-handle along the dotted line indicated in the figure).  \qed
  
  \vspace{.1in} 
  
    Next we check whether  $R$ and $S: \partial W(0,1) \to \partial W(0,1)$ are cork automorphisms?.  In previous cork examples,  this was done by use of the adjunction inequality (e.g.\cite{a3}). More specifically, finding a slice knot $\gamma \subset W(0,1)$ so that the image of $\gamma$ under $R$ can not be slice, would imply $R$ is a cork automorphism. Recently, in \cite{dhm}, by using another technique, it was showned that these maps were cork automorphisms. So we won't pursue the adjunction technique here, instead we will answer the question of whether the homotopy sphere obtained by doubling $W(0,1)$  by the maps $R$ or $S$ is standard $S^{4}$? 
    
 \vspace{.1in}   
    
\begin{figure}[ht]  \begin{center}
 \includegraphics[width=.69\textwidth]{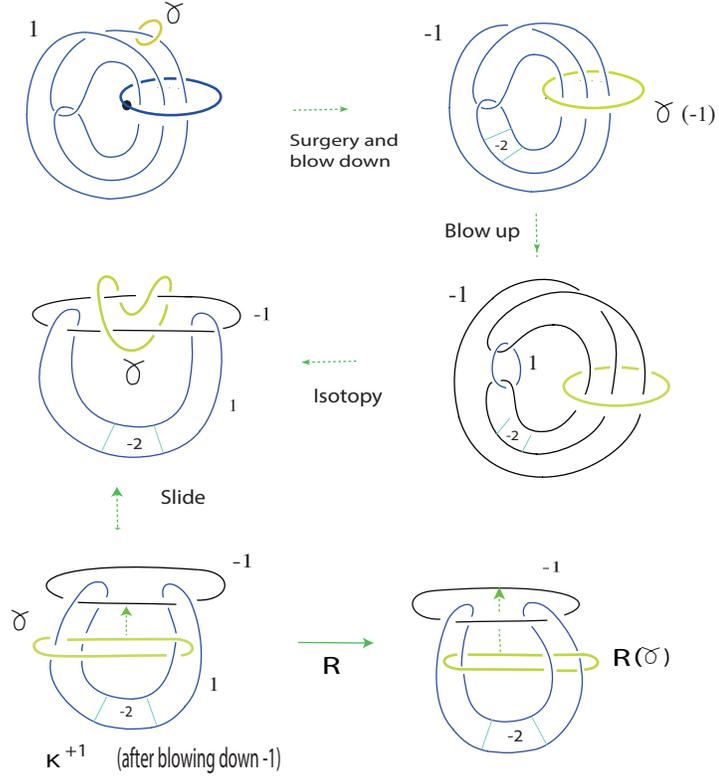} 
    \vspace{.05in}      
\caption{Tracing dual circle of $2$-handle of $W(0,1)$ to $K^{+1}$}      \label{c4} 
\end{center}
  \end{figure}

\begin{figure}[ht]  \begin{center}
   \vspace{-.1in}
 \includegraphics[width=.65\textwidth]{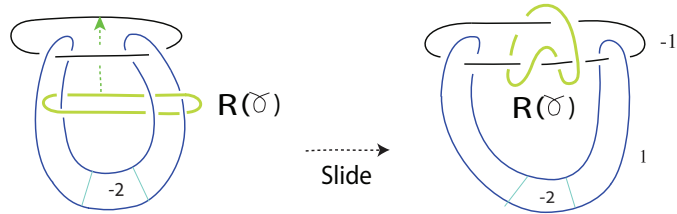}       
\caption{Sliding $R(\gamma)$ over the $2$-handle}      \label{c46} 
\end{center}
  \end{figure}

\begin{figure}[ht]  \begin{center}
 \includegraphics[width=.7\textwidth]{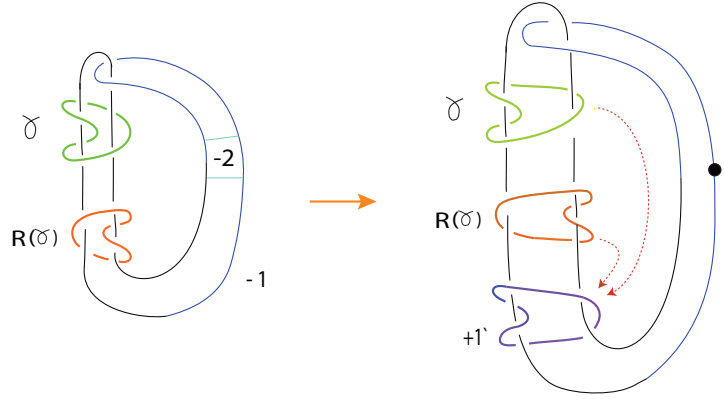}       
\caption{Positions of $\gamma$ and $R(\gamma)$ in $\partial W(0,1)$}      \label{c47} 
\end{center}
  \end{figure}

\begin{figure}[ht]  \begin{center}
 \includegraphics[width=.7\textwidth]{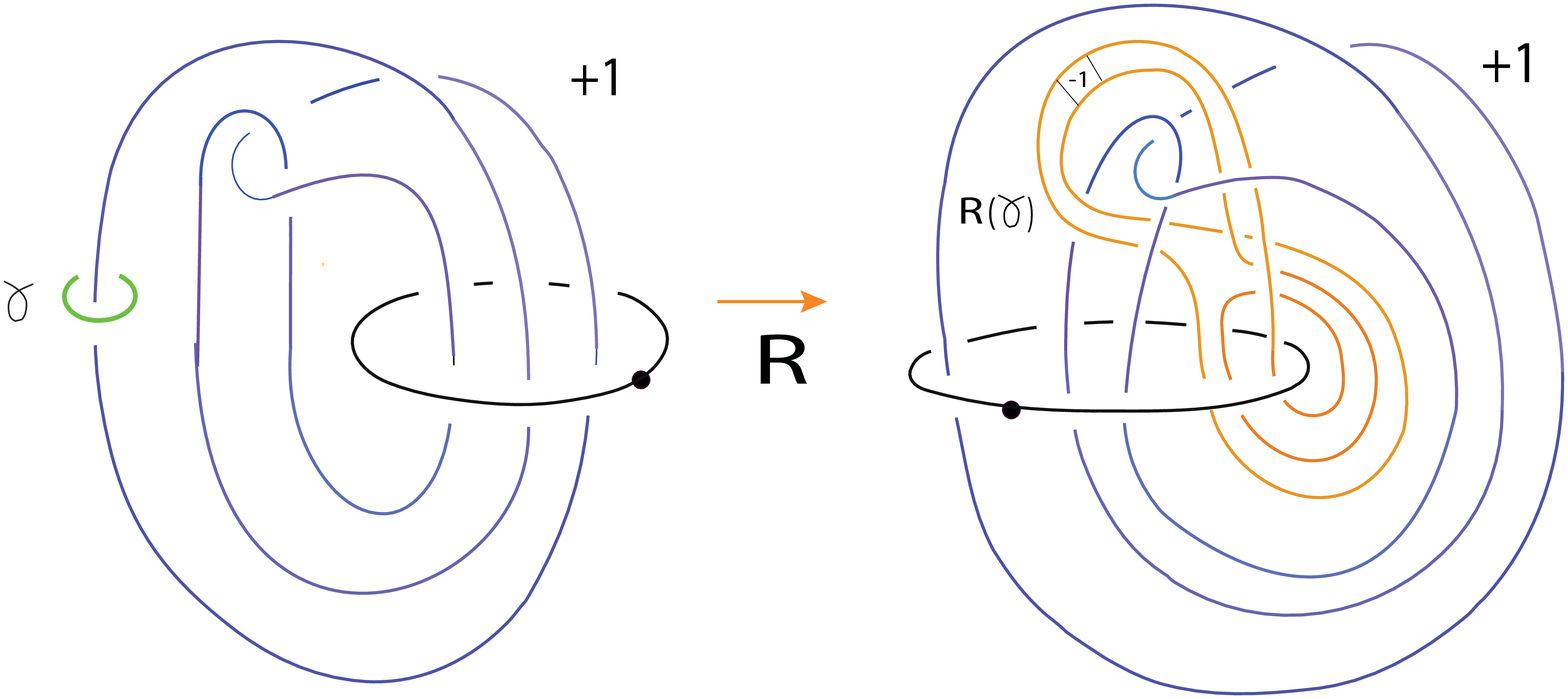}       
\caption{$R:\partial W(0,1)\to \partial W(0,1)$}     \label{c50} 
\end{center}
  \end{figure}

    Amazingly, the solution of this goes through the same interesting steps as in the solution of the Cappell-Shaneson homotopy sphere problem, which took 30 years to settle \cite{a4}.    
   For this we trace image  $R(\gamma)$ of the dual circle 
  $\gamma$ of the $2$-handle of $W(0,1)$ to  $\partial K^{+1}$. Also it follows from the construction that $R(\gamma)= S(\gamma)$.  Figure~\ref{c4} shows how to identify $\partial W(0,1) $ with $K^{+1}$ (just blow-down $-1$ framed unknot in the picture). This figure  identifies the action of $R$ as $180^{o}$ rotation. After sliding $R(\gamma)$ over the $2$-handle we get Figure~\ref{c46}. In Figure~\ref{c47} we draw the positions of $\gamma$  and $R(\gamma)$   in the same picture of $W(0,1)$.  Then after sliding $\gamma$  and $R(\gamma)$ over the $2$-handle (along the indicated arrows) we arrive Figure~\ref{c50}, which describes the action of R on $\partial W(0,1)$, in particular  indicating the position of $R(\gamma)$. 
    
    
  \clearpage
  
   \subsection{Proof of Theorem~\ref{hs}}
   
\begin{figure}[ht]  \begin{center}
 \includegraphics[width=.9\textwidth]{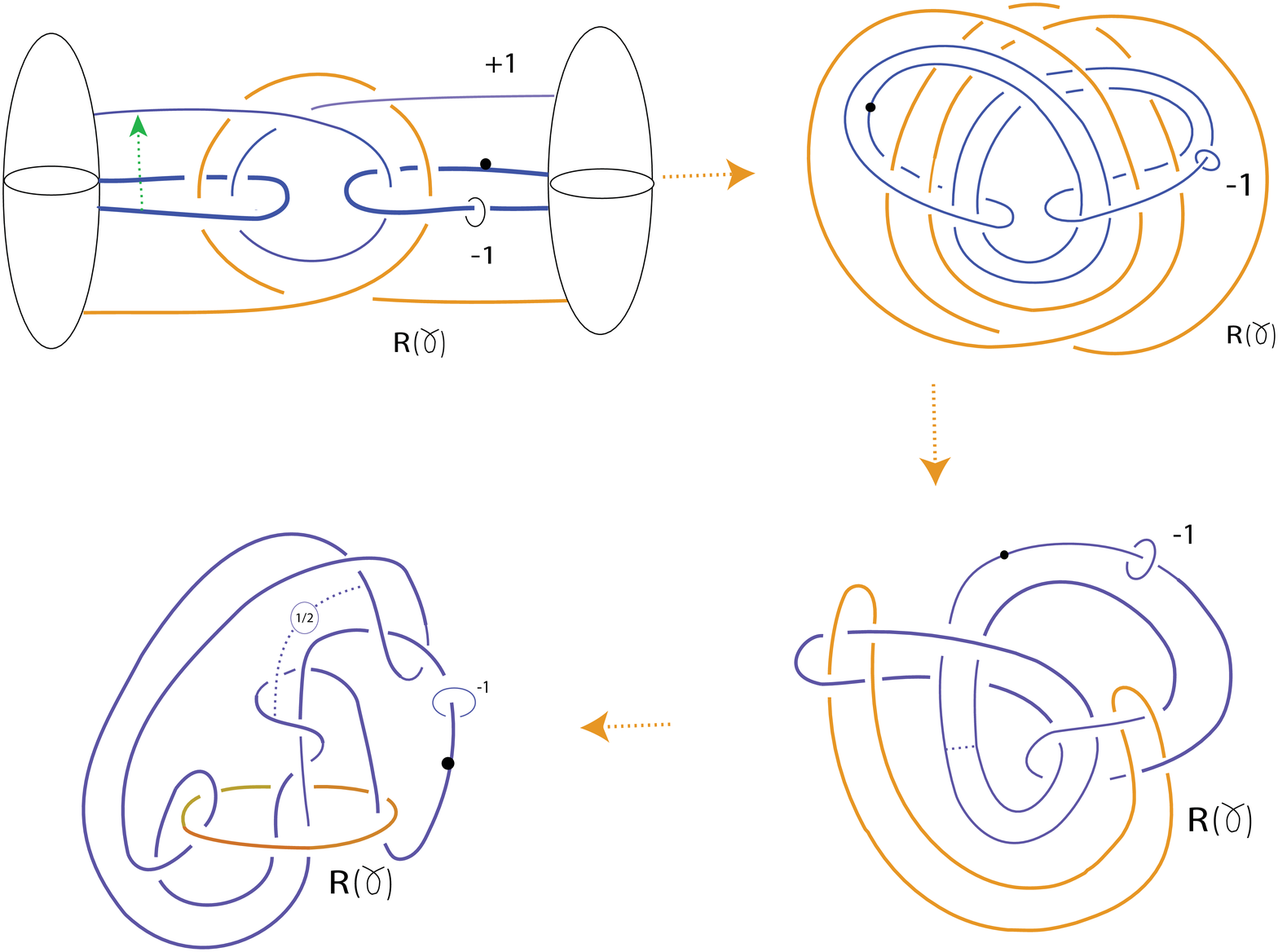}       
\caption{Tracing $R(\gamma)\subset \partial W(0,1)$}      \label{c51} 
\end{center}
  \end{figure}

   Clearly, up to $3$-handles, the handlebody of $\Sigma =- W(0,1)\smile _{R(\gamma)}W(0,1)$ can be obtained by attaching $0$-framed $2$-handle to $W(0,1)$ along $R(\gamma)$, which is the dual  $2$-handle of $-W(0,1)$. This means attaching a $2$-handle to the second picture of Figure~\ref{c50} along $R(\gamma)$.  Careful reader will notice that the big left twist between the curves in the middle of Figure~\ref{c50} can represented in a simpler way by introducing a cancelling $1/2$- handle pair (a dotted circle and small linking  $-1$ framed circle). This gives the first picture of Figure~\ref{c51}. 
   

\begin{figure}[ht]  \begin{center}
 \includegraphics[width=.7\textwidth]{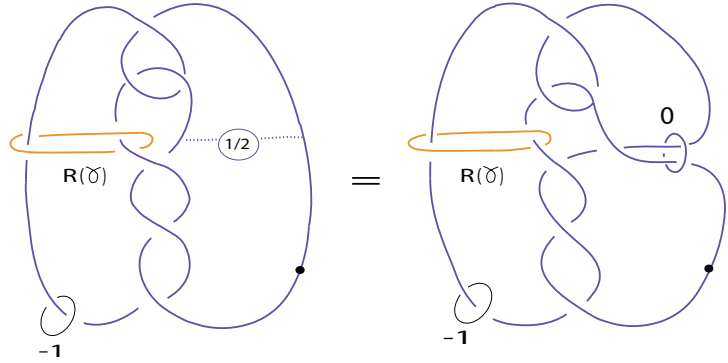}       
\caption{Two different ribbons describing $S^{2}\subset S^{4}$}      \label{c52} 
\end{center}
  \end{figure}

   This observation will now lead us to Gluck construction along an imbedded $2$-sphere $S^{2}\subset S^{4}$: Now we slide the middle $1$-handle over the $2$-handle as indicated in this figure (note that $1$-handles can not slide over the $2$-handles, unless they are in the form of Figure 1.16 of \cite{a3}). Then proceed to the last picture of Figure~\ref{c51}. The dotted arcs in Figure~\ref{c51} are the ribbon moves reminding us the bounding disk of the ribbon $1$-handle. Drawing the last picture of Figure~\ref{c51} in a more symmetric way we get the first picture of Figure~\ref{c52}.

   This last picture reveals an amazing fact: It describes a Gluck twisted $S^4$ twisted along an imbedded $S^{2}\subset S^{4}$, where this $2$-sphere is obtained by putting together the two different ribbon disks $D_{\pm}\subset B^{4}_{\pm}$, which the Stevedore knot $K$ bounds. The two distinct ribbon moves are related to each other by $180^{o}$ rotation of $S^{3}$ as indicated by the picture. Second picture of Figure~\ref{c52} (after ribbon move performed) is the picture of the Glucked $2$-sphere $S^2\subset S^4$. $R(\gamma)$ represents a zero framed $2$-handle.

          \begin{figure}[ht]  \begin{center}
 \includegraphics[width=.4\textwidth]{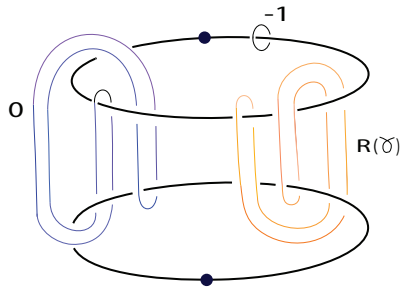}       
\caption{Another view of $R(\gamma)\subset W(0,1)$. Attaching $2$-handle to $R(\gamma)$ turns picture to Gluck twisted $S^{4}$}      \label{c53} 
\end{center}
  \end{figure}

\begin{figure}[ht]  \begin{center}
 \includegraphics[width=.9\textwidth]{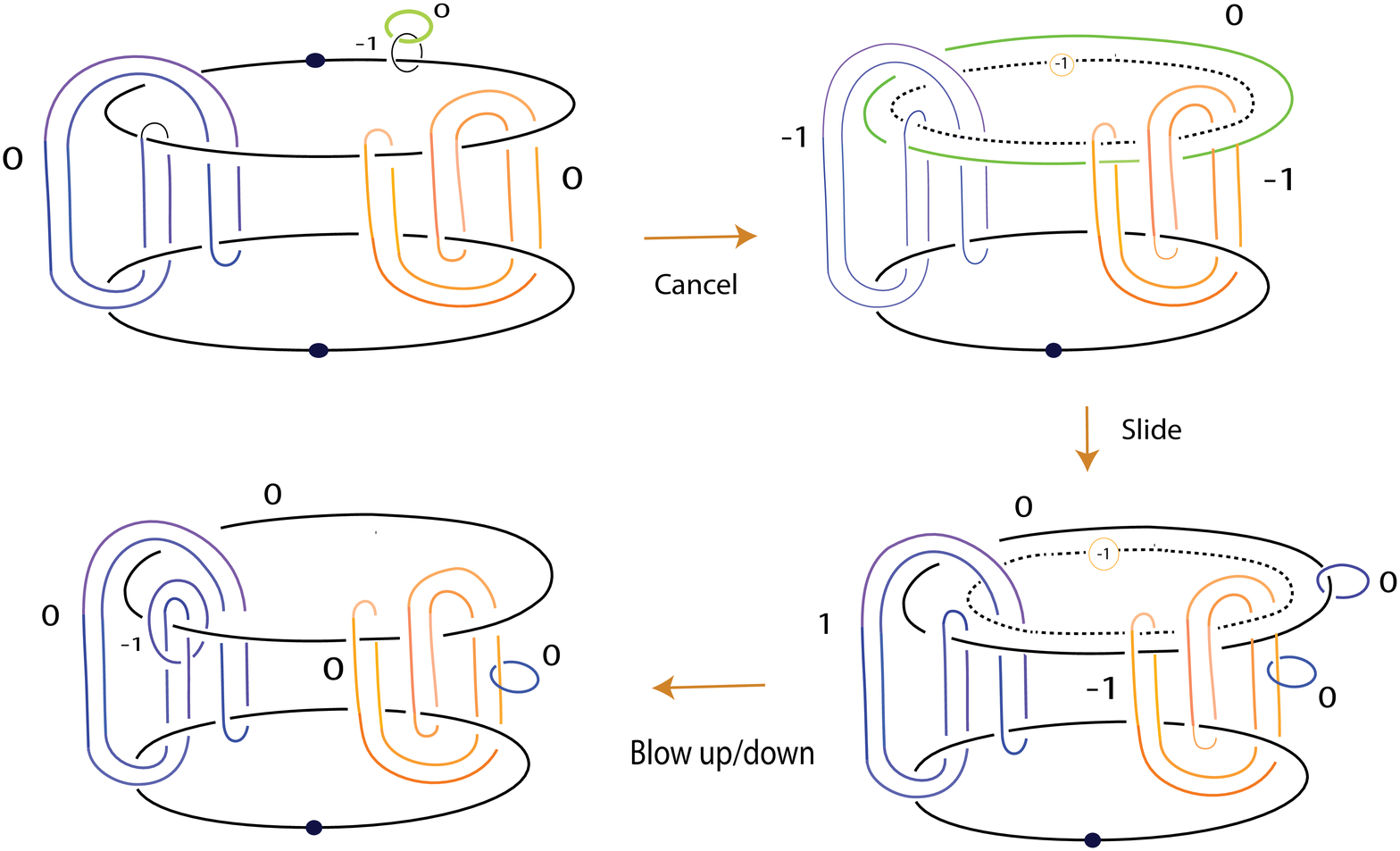}       
\caption{Using Gluck twist to cancel $1$-handle}      \label{c54} 
\end{center}
  \end{figure}

  
 Figure~\ref{c53} is the same as the second picture of Figure~\ref{c52}, drawn more symmetric way. It is easy to check that in this picture the linking circle of the small $-1$ framed circle is unknot in $S^{3}$ boundary.  Now by using this unknot, we can attach a $2/3$ - canceling handle pair (new $2$-handle is the $0$-framed small red circle in the first picture of Figure~\ref{c54}).

 Next we use the trick from (Figure 14.11 of \cite{a3}):  That is, we cancel the circle-with-dot at top with its linking $-1$ framed circle and get the second picture of Figure~\ref{c54}. Then after the obvious handle slide over the large $0$-framed $2$-handle at the top right picture of Figure~\ref{c54}, we obtain the third picture of Figure~\ref{c54},  where we can now see a cancelling $1/2$-handle pair. So this picture can be thought of a handlebody without $1$-handles (i.e. it consists of two $2$-handles and two $3$-handles), and hence turning it upside down we will give get a handlebody without $3$-handles!  Having noted this, we can now turn this handlebody upside down (as the process described in \cite{a3}).

       \begin{figure}[ht]  \begin{center}
 \includegraphics[width=.74\textwidth]{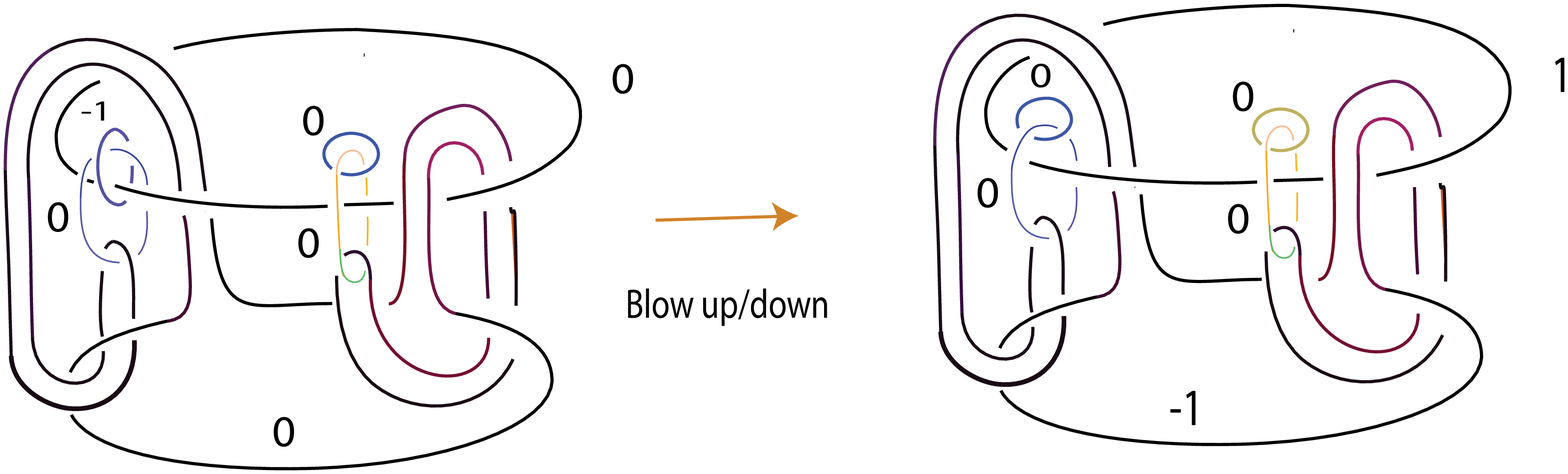}       
\caption{Starting the process of turning upside down}      \label{c55} 
\end{center}
  \end{figure}

  \newpage
  
\begin{figure}[ht]  \begin{center}
 \includegraphics[width=.7\textwidth]{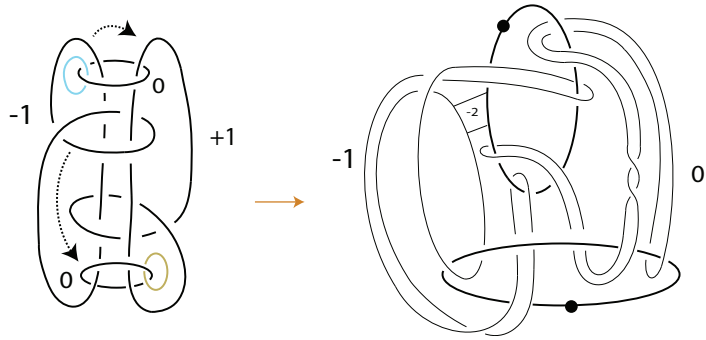}       
\caption{Turning upside down process}      \label{c56} 
\end{center}
  \end{figure}

\begin{figure}[ht]  \begin{center}
 \includegraphics[width=.85\textwidth]{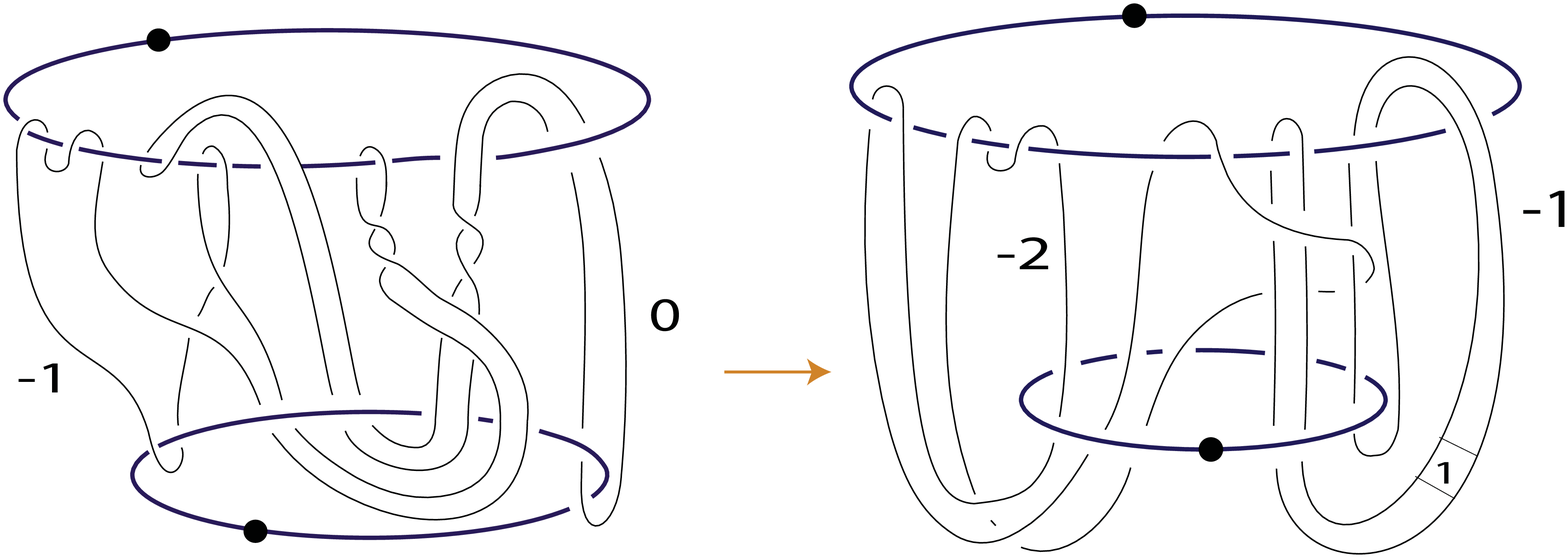}       
\caption{$\Sigma $ without $3$-handles }      \label{c57} 
\end{center}
  \end{figure}

\begin{figure}[ht]  \begin{center}
 \includegraphics[width=.75\textwidth]{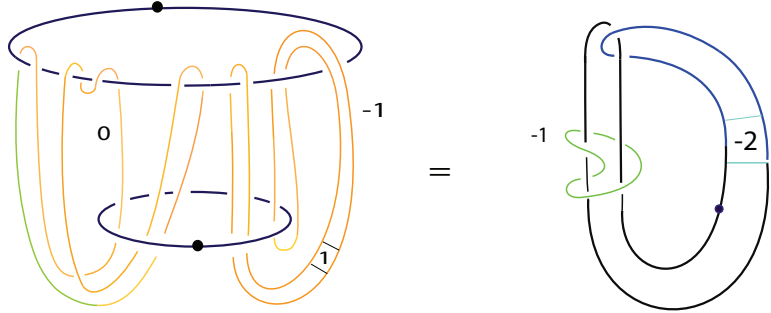}       
\caption{two equivalent handlebody descriptions of $\Sigma$}      \label{c58} 
\end{center}
  \end{figure}

\begin{figure}[ht]  \begin{center}
 \includegraphics[width=.85\textwidth]{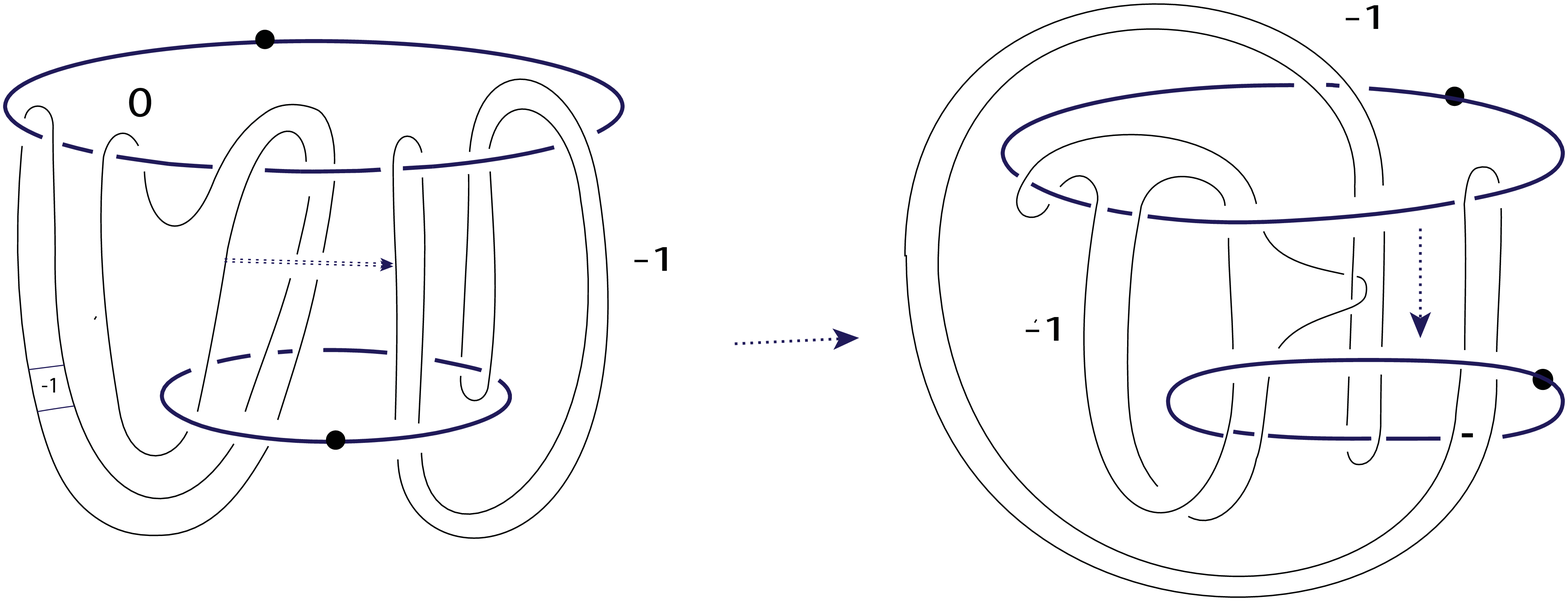}       
\caption{Sliding handles to simpilfy $\Sigma$}      \label{s4} 
\end{center}
  \end{figure}

\begin{figure}[ht]  \begin{center}
 \includegraphics[width=.85\textwidth]{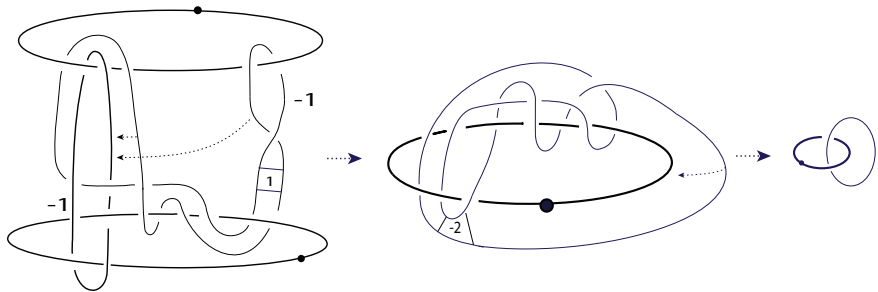}       
\caption{More handle slides to identify $\Sigma \approx S^4$}      \label{s5} 
\end{center}
  \end{figure}

  For this, we ignore the cancelling $1/2$ handle pair, and carry the duals of the remaing $2$-handles to the boundary of $\# 2(S^{1}\times B^{3})$ by any diffeomorphism. Starting with the third picture of Figure~\ref{c54}, we proceed from Figures~\ref{c55} to Figure~\ref{c57} by self described steps similar to in Figure 14.11 of \cite{a3}, and end up with Figure~\ref{c57}, and the first picture of Figure~\ref{c58}.  By canceling a $1/2$-handle pair gives the second picture of Figure~\ref{c58} (dotted Stevedore knot represents a ribbon $1$-handle).  Finally, by performing the indicated handle slides to Figure~\ref{s4} we obtain Figure~\ref{s5}, and from Figure~\ref{s5} we arrive to $S^4$ ($4$-handle is not drawn). \qed

\clearpage

  {\it Acknowledgement:}   
We thank Danny Ruberman for helpful remarks.


\begin{thebibliography}{99999}

  \bibitem[A1]{a1} S. Akbulut, {\em A Fake compact contractible 4-manifold}, J. Differ. Geom. 33 (1991) 335--356.
 
  \bibitem[A2]{a2}  S. Akbulut, {\em On infinite order
  corks},  Proc. {G}\"{o}kova {G}eometry-{T}opology {C}onf.,  (2017) 151--157.
   
 \bibitem [A3]{a3} S. Akbulut, {\em $4$-Manifolds}, vol.~25
  of Oxford Graduate Texts in Mathematics, Oxford University Press, Oxford, 2016.
  
   \bibitem [A4]{a4} S. Akbulut, {\em Cappell-Shaneson homotopy spheres are standard}, Ann. of Math, 171 (2010) 2171-2175.
  
  \bibitem [A5]{a5} S. Akbulut, {\em Homotopy $4$-spheres associated to an infinite order loose cork},
      https://arxiv.org/pdf/1901.08299.pdf
       
 \bibitem [A6]{a6} S. Akbulut, {\em An exotic $4$-manifold}, 
Journ. of Diff. Geom. 33, (1991), 357-361.

   \bibitem [AK]{ak} S. Akbulut and R. Kirby, {\em  Mazur manifolds}, Mich. Math. J. 26 (1979) 259--284.
   
 \bibitem[AM]{am} S. Akbulut and R. Matveyev, {\em A convex decomposition theorem for 4-manifolds}, Int. Math. Res. Notices, no. 7 (1998) 371-381.
 
 \bibitem[AY]{ay} S. Akbulut and K. Yasui, {\em Corks, Plugs and exotic structures}, Journal of G\"{o}kova Geometry Topology, volume \textbf{2} (2008), 40--82.

  \bibitem[BO]{bo} M. Boileau and J. P. Otal, {\em Scindements de Heegaard et groupe de hom\'{e}otopies
des petites vari\'{e}t\'{e}s de Seifert}, Invent. Math. {\bf 106},  (1991), 85--108.
  
  \bibitem[BW]{bw} M. Brittenham and Y.-Q. Wu {\em The classification of exceptional Dehn surgeries on 2-bridge knots},  Comm. Anal. Geom. {\bf 9} (2001), no. 1, 97--113.
  
   \bibitem[CFHS]{cfhs} C.L. Curtis, M.H. Freedman, W.C. Hsiang, R. Stong, {\em A decomposition theorem for $h$-cobordant smooth simply-connected compact $4$-manifolds}, Invent. Math. {\bf 123}, no.2 (1996), 343--348.  
      
 \bibitem[DHM]{dhm} I. Dai, M. Hedden, A. Mallick {\em Corks, involutions, and Heegaard Floer homology}, \\ arXiv:2002.02326 

 \bibitem[G]{g} R. E. Gompf, {\em Infinite order corks},  Geom. Topol. {\bf 21} (2017), no. 4, 2475--2484.

  \bibitem[M]{m} R. Matveyev, {\em A decomposition of smooth simply-connected $h$-cobordant $4$-manifolds}, Jour. Diff. Geom. {\bf 44},  no. 3 (1996), 571--582.
   
    \bibitem[MT]{mt} T.E. Mark, B. Tosun,  {\em Obstructing pseudoconvex embeddings and contractible Stein fillings for Brieskorn spheres}, Adv. Math. {\bf 335} (2018), 878--895.
      
    \bibitem[R]{r} D. Ruberman, {\em Symmetries of graphs and automorphisms of 3-manifolds}. 
 
  \end{thebibliography}
\end{document}